\documentclass[reqno]{amsart}
\usepackage{txfonts}
\usepackage{bbm}
\usepackage{mathrsfs}
\usepackage{amsmath}
\usepackage{paralist}
\usepackage{graphics}
\usepackage{epsfig}
\usepackage[pdftex]{hyperref}
\usepackage{amsfonts}
\usepackage{CJK}
\usepackage{fancyhdr}
\usepackage{graphicx}
\usepackage[dvips,usenames]{color}
\usepackage{titletoc}
\usepackage{latexsym}
\usepackage{amssymb}
\usepackage{multicol}
\usepackage{graphics}
\usepackage{subfigure}
\usepackage{indentfirst}
\usepackage{cases}
\usepackage{curves}
\usepackage{cite}
\newtheorem{theorem}{Theorem}[section]

\newtheorem{lemma}{Lemma}[section]
\newtheorem{definition}{Definition}[section]
\newtheorem{proposition}{Proposition}[section]
\newtheorem{remark}{Remark}[section]
\newtheorem{example}{Example}[section]

\numberwithin{equation}{section}

\def\x#1{(\ref{#1})}
\def\R{{\Bbb R}}

\def\N{{\Bbb N}}

\def\bc{\begin{center}}
\def\ec{\end{center}}
\def\ba{\begin{array}}
\def\ea{\end{array}}
\def\be{\begin{equation}}
\def\ee{\end{equation}}
\def\bea{\begin{eqnarray}}
\def\eea{\end{eqnarray}}
\def\beaa{\begin{eqnarray*}}
\def\eeaa{\end{eqnarray*}}

\def\ben{\begin{enumerate}}
\def\een{\end{enumerate}}
\def\hh{\!\!\!\!}

\def\EQ{\hh & = & \hh}

\def\LE{\hh & \le & \hh}
\def\GE{\hh & \ge & \hh}

\def\oo{\infty}
\def\ifl{\iffalse}
\def\lb{\label}
\def\prf{\mbox{\bf Proof.~}}

\title[]{Nonuniform   Dichotomy  Spectrum Intervals: Theorem and Computation}

\author[H. Zhu]
{Hailong Zhu $^{1}$}

\address{$^1$ School of Statistics and Applied Mathematics, Anhui
University of Finance and Economics, Bengbu 233030, China}

\email{hai-long-zhu@163.com (H. Zhu)}


\begin{document}

\begin{abstract} Under the condition of nonuniformly bounded growth, 
the relationship of the nonuniform exponential dichotomy spectrum and the other two classical  spectrums (the Lyapunov spectrum and Sacker-Sell spectrum) is given, and the stability of these spectrums under small linear perturbations are summarized and presented in this paper. A main goal of this paper is to discuss the theory for the computation of these spectrums under the condition of nonuniformly bounded growth, and this extends the work of Dieci and Vleck \cite{dv-02}, which compute the Lyapunov spectrum and Sacker-Sell spectrum under the condition of bounded. Finally, an example is given to illustrate and verify the theoretical results.
\end{abstract}

\maketitle

\section{\bf{Introduction}}
\setcounter{equation}{0} \noindent

Lyapunov exponents was introduced by Lyapunov himself, reprinted in \cite{lv-92}.
In this paper, Lyapunov exponents was generalized  for illustrating the characterization of exponential growth rates of  time varying matrix functions. For an $n$-dimensional (time varying) differential equations,  there are $n$ Lyapunov exponents, and it is natural to think about

Since then, different characterizations of spectrums for  linear nonautonomous differential equation have been proposed. Among them, one of the most famous spectrums is
dichotomy spectrum (or called Sacker-Sell spectrum, dynamical spectrum), which was introduced by Sacker and Sell in \cite{ss1, ss2} defined by exponential dichotomies to study the linear skew product flows. Since these classical works,  a lot of research has been done to understand and extend this fruitful concept in various ways. For
example, A spectral  theory about linear difference equations has been studied in \cite{amz, am, as01, as02}. Reducibility and normal forms  for nonautonomous
differential equations by using dichotomy spectrum has been given in \cite{sie2-02,sie3-02}. For more results about dichotomy spectrum,  see \cite{dp-11,dk-12,pot09,pot12, sie-02}  and the references therein.

In the computation of spectral intervals, both for dichotomy spectrum $\Sigma_{ED}(A)$ and Lyapunov spectrum $\Sigma_{L}(A)$, SVD and QR methods have been proposed by Dieci and Vleck \cite{dv-02, dv2-02} to study the computation methods for these spectrums. After that, further research on this topic has been proposed by Dieci and his collaborators (see \cite{d-06, dv2-07, d-08, d-11} for details). For more information about the theoretical and numerical analysis of dichotomy spectrum, one can refer to \cite{f-13,h-10} and the references cited therein.

On  the   other  hand,  as  Barreira  and   Valls  mentioned  in \cite{bv2-08},  the   classical notion   of   exponential  dichotomy  substantially   restricts   some  dynamics,   and from  the   point  of   view   of   ergodic  theory,  almost   all  linear   variational   equations  have   a nonuniform  exponential  behavior. During the  last  several  years,
 a more generalized  concept: nonuniform  exponential  dichotomy has been introduced and investigated by Barreira  and   Valls (see e.g., \cite{bv-06,bv-07,bv-08}).  Based  on  the   study   of exponential  dichotomy,  the  nonuniform dichotomy  spectral   theory
was  introduced in \cite{chu-15,zh-14} for linear nonautonomous system with the coefficients being nonuniformly bounded growth (see Definition \ref{def41} below).

Here we mention that  the numerical methods proposed by Dieci and Vleck in \cite{dv-02, dv2-02} demands  the coefficients of the linear systems to be bounded.  Otherwise, the numerical technique for computing Sacker-Sell spectrum, which is based on the condition of integral separateness, is not quite right.  For example, consider  the following two dimensional diagonal system
\be\lb{b12} \left(\begin{array}{c}
         \dot{x}_{1} \\
         \dot{x}_{2}
       \end{array}\right)=\left(
  \begin{array}{cc}
    \omega_{1} & 0\\
    0 & \omega_{2}t\sin t \\
  \end{array}
\right)\left(\begin{array}{c}
         x_{1} \\
         x_{2}
       \end{array}\right)
\ee
with $\omega_{1}>\omega_{2}>0$ be real paraments. One can see that the coefficients of \x{b12} is not bounded, and \x{b12} is not integrally separated (see \cite{zc}). Moreover, one can prove that the dichotomy spectrum $\Sigma_{ED}$ of \x{b12}  is trivial, i.e., $\Sigma_{ED}=\R$, and the nonuniformly  dichotomy spectrum is $\Sigma_{NED}=\{\omega_{1}\}\cup [-\omega_{2}, \omega_{2}]$ (see Example 2.1 in \cite{chu-15} for details, Remark \ref{remfj1} below also presents an explanation from the point of view of numerical analysis).

This work, inspired by both the classical notion of  dichotomy spectrum \cite{ss1,ss2} and the notion of nonuniform  dichotomy  spectrum introduced by \cite{chu-15,zh-14},  is an attempt to discuss the relationship of three different spectrums:  $\Sigma_{L}(A)$, $\Sigma_{ED}(A)$ and $\Sigma_{NED}(A)$, and extend the numerical technique  developed by Dieci and Vleck \cite{dv-02, dv2-02}  for studying  linear nonautonomous system with the coefficients being nonuniformly bounded growth.

An outline of the paper is as follows. In Section 2, the basic definitions and properties of $\Sigma_{L}(A)$, $\Sigma_{ED}(A)$ and $\Sigma_{NED}(A)$ will be presented. Section 3 discusses the relationship of $\Sigma_{L}(A)$, $\Sigma_{ED}(A)$ and $\Sigma_{NED}(A)$, and summarizes the stability of these spectrums under small linear perturbations.
In Section 4, we first establish necessary and sufficient condition of Steklov function and weak integral separateness under the condition of nonuniformly bounded growth. Thus we can use this relationship to show the numerical methods for  $\Sigma_{ED}(A)$ and $\Sigma_{NED}(A)$. An example will be given in Section 5 to illustrate and verify the theoretical results.

\section{\bf{Lyapunov, exponential dichotomy, and nonuniform exponential dichotomy spectrum}}
\setcounter{equation}{0} \noindent
Given an $n$-dimensional linear system
\be\lb{b1}
\dot{x}=A(t)x,
\ee
where $x(t)\in \R^{n}$ and $A(t)$: is a $n\times n$ matrix with real entries depending
continuously on $t \ge 0$. Consider  the trivial solution of \x{b1}. It is well known that if the matrix function
$A(t)$ is constant, i.e., $A(t) = A$ for all $t \ge 0$, then the zero solution of \x{b1} is asymptotically (and indeed, exponentially) stable if and only if the real part
of every eigenvalue of the matrix $A$ is negative. A similar result holds in the
case when the matrix function $A(t)$ is periodic by using the Floquet theory. For the general (nonautonomous) case, we need to consider the spectral intervals instead of eigenvalues, so in this section we first recall the definitions of the next two classical concepts of spectrum: Lyapunov spectrum  $\Sigma_{L}(A)$, exponential dichotomy spectrum  $\Sigma_{ED}(A)$, and then we introduce a third related one, the nonuniform exponetial dichotomy spectrum, $\Sigma_{NED}(A)$.

 \vspace{6pt}
 \subsection{Lyapunov spectrum.} Given a fundamental matrix solution $\Phi(t)$ of \x{b1}, define $\lambda_{j}, j=1,\ldots,n$, as
 \[
 \lambda_{j}(\Phi(t)):=\limsup_{t\rightarrow +\oo}\frac{1}{t}\ln\|\Phi(t)e_{j}\|,
 \]
where the vector norm is the $2$-norm (invariant under orthogonal
transformations), and the $e_{j}$ is the unit column-vector in the $x_{j}$ direction, i.e., \[e_{j}=(\underbrace{0,\ldots,0,1}_{j},0,\ldots,0)^{T}.\]
When the sum of the numbers $\lambda_{j}$ is minimized as we vary over all possible fundamental matrix solutions of the system, i.e.,
\[\sum_{j=1}^{n}\lambda_{j}^{s}:=\inf\left\{\sum_{j=1}^{n}\lambda_{j}(\Phi(t)):~\Phi(t)~is~a~fundamental~matrix~of ~\x{b1}\right\}\]
 then the numbers $\lambda_{j}^{s}, j=1,\ldots,n$ are called \emph{(upper) Lyapunov exponents}, and the corresponding basis is called \emph{normal}.

 Now we consider the linear differential equation which is dual to  \x{b1}
 \be\lb{b3}\dot{y}=-A^{T}(t)y.\ee
 Similarly, one can have the upper Lyapunov exponents $\lambda_{j}^{i}, j=1,\ldots,n$ of \x{b3}, which are the lower Lyapunov exponents of \x{b1} (see e.g., \cite{dv-02} for details). Let $\lambda_{j}^{s}, \lambda_{j}^{i}$ be ordered: $\lambda_{1}^{s}\ge\lambda_{2}^{s}\ge \cdots \ge \lambda_{n}^{s}$ and $\lambda_{1}^{i}\ge\lambda_{2}^{i}\ge \cdots \ge \lambda_{n}^{i}$.
 Considering such a fact that $\lambda_{j}^{i}\le \lambda_{j}^{s}$, then the Lyapunov spectrum can be defined as
 \[
 \Sigma_{L}:=\bigcup_{j=1}^{n}[\lambda_{j}^{i},\lambda_{j}^{s}].\]
Especially, the system is called \emph{regular} while $\lambda_{j}^{i}=\lambda_{j}^{s}$ for all $j=1,\ldots, n$.

 \vspace{6pt}
 \subsection{Exponential dichotomy spectrum.}
Recall that \x{b1} admits   an \emph{exponential  dichotomy} if there exist an invariant projection  $P$, and constants $\alpha>0$, $M>0$ such that
 \begin{equation*}
\|\Phi(t)P\Phi^{-1}(s)\|\le M e^{-\alpha(t-s)}, \quad {\rm for} \quad 0\le s \le t,
\end{equation*}
and
\begin{equation*}
\|\Phi(t)Q\Phi^{-1}(s)\|\le M e^{\alpha(t-s)}, \quad {\rm for} \quad 0\le t \le s,
\end{equation*}
where  $Q =I_{n}-P$ is  the   complementary   projection. Furthermore, for any fixed $\gamma\in \R$, write a   shifted   system
\begin{equation}\lb{b7}
  \dot{x}=\left[A(t)x-\gamma I_{n}\right]x.
\end{equation}
Then the   \emph{exponential dichotomy  spectrum}  of {\rm \x{b1}} is given by  the   set
\[\Sigma_{NED}(A)=\{\gamma\in \R: ~{\rm \x{b7}}~admits~no ~ exponential~ dichotomy\},\]
and   the   resolvent   set $\rho_{ED}(A)=\R \setminus \Sigma_{ED}(A)$ is its complements.

 In \cite{ss2, sie-02}, it has been shown that $\Sigma_{ED}(A)$ is at most a disjoint union of $n$ closed intervals. This means that $\Sigma_{ED}(A)=\emptyset$ or $\Sigma_{ED}(A)=\R$ or $\Sigma_{ED}(A)$ is   in   one   of   the  four   cases
\[\Sigma_{ED}(A)=\left\{\begin{array}{c}
                   [a_{1},b_{1}] \\or\\
                   (- \oo,b_{1}]
                 \end{array}\right\}\cup [a_{2},b_{2}]\cup \cdots \cup [a_{k-1},b_{k-1}]\cup
\left\{ \begin{array}{c}
  [a_{k},b_{k}] \\
  or
  \\

  [a_{k},  \oo)
\end{array} \right\}
\]
for some $k: 1\le k\le n$.

\vspace{6pt}
 \subsection{Nonuniform exponential dichotomy spectrum.} In \cite{bv-06, bv-07}, Barreira and Valls propose a new notion
called nonuniform, which extends the notion of dichotomy of uniform. Later, \cite{chu-15} presents a new spectrum for \x{b1} based upon the nonuniform exponential dichotomy.

Recall that \x{b1}  admits   a   \emph{nonuniform  exponential  dichotomy} if there exist an invariant projection  $P$, constants $\alpha>0$, $M>0$, and $\varepsilon \in [0,\alpha)$ such that
 \begin{equation}\lb{b8}
\|\Phi(t)P\Phi^{-1}(s)\|\le M e^{-\alpha(t-s)}e^{\varepsilon s}, \quad {\rm for} \quad 0\le s \le t,
\end{equation}
and
\begin{equation}\lb{b9}
\|\Phi(t)Q\Phi^{-1}(s)\|\le M e^{\alpha(t-s)}e^{\varepsilon s}, \quad {\rm for} \quad 0\le t \le s,
\end{equation}
where  $Q =I_{n}-P$ is  the   complementary   projection.
Then the   \emph{nonuniform exponential  dichotomy  spectrum}  of {\rm \x{b1}} is  given by the   set
\[\Sigma_{NED}(A)=\{\gamma\in \R: ~{\rm \x{b7}}~admits~no~nonuniform~ exponential~ dichotomy\},\]
and   the   resolvent   set $\rho_{NED}(A)=\R \setminus \Sigma_{NED}(A)$ is its complements.

Similarly, it has been shown in \cite{chu-15} that $\Sigma_{NED}(A)$ is at most a disjoint union of $n$ closed intervals. This means that $\Sigma_{NED}(A)=\emptyset$ or $\Sigma_{NED}(A)=\R$ or $\Sigma_{NED}(A)$ is   in   one   of   the  four   cases
\[\Sigma_{NED}(A)=\left\{\begin{array}{c}
                   [a_{1},b_{1}] \\or\\
                   (- \oo,b_{1}]
                 \end{array}\right\}\cup [a_{2},b_{2}]\cup \cdots \cup [a_{k-1},b_{k-1}]\cup
\left\{ \begin{array}{c}
  [a_{k},b_{k}] \\
  or
  \\

  [a_{k},  \oo)
\end{array} \right\}
\]
for some $k: 1\le k\le n$.

\section{\bf{Relationship of spectrums $\Sigma_{L}(A)$, $\Sigma_{ED}(A)$ and $\Sigma_{NED}(A)$.}}
\setcounter{equation}{0} \noindent
 It is well known that the notion of  Lyapunov exponents, exponential dichotomy together with some of their variants, extensions, and modifications, play a central role in the
study of general theory of dynamical systems.
 To gain insight into the behavior of the dynamical approaches of \x{b1}, several aspects are discussed in this section to illustrate the relationship of spectrums $\Sigma_{L}(A)$, $\Sigma_{ED}(A)$ and $\Sigma_{NED}(A)$.

 We first present the relation of inclusion of these three spectrums.
 \begin{proposition}
   For an $n$-dimensional linear system {\rm\x{b1}}, we have the following chain of implications
   \[\Sigma_{L}(A)\subset\Sigma_{NED}(A)\subset\Sigma_{ED}(A).\]
 \end{proposition}
 \prf{Clearly $\Sigma_{NED}(A)\subset\Sigma_{ED}(A)$ due to the fact $\varepsilon\ge 0$ in \x{b8}-\x{b9} (see \cite{chu-15} for details). Now we prove that $\Sigma_{L}(A)\subset\Sigma_{NED}(A)$. Obviously,  $\Sigma_{L}(A)\subset\Sigma_{NED}(A)$ if $\Sigma_{NED}(A)=\R$. Conversely, if $\Sigma_{NED}(A)=\emptyset$, then $\rho_{NED}(A)=\R$. This means that for any $\lambda\in \R$, there exist an invariant projection  $P$, constants $\alpha>0$, $M>0$, and $\varepsilon \in [0,\alpha)$ such that
 \[
\|\Phi_{\lambda}(t)P\Phi_{\lambda}^{-1}(s)\|\le M e^{-\alpha(t-s)}e^{\varepsilon s}, \quad {\rm for} \quad 0\le s \le t,
\]
or  equivalently,
\[
\|\Phi(t)P\Phi^{-1}(s)\|\le M e^{(\lambda-\alpha)(t-s)}e^{\varepsilon s}, \quad {\rm for} \quad 0\le s \le t.
\]
Set $s=0$, the inequality above implies that
\[\|\Phi(t)P\Phi^{-1}(0)x_{0}\|\le M \|x_{0}\|  e^{(\lambda-\alpha)t} \quad {\rm for} \quad 0 \le t.
\] with any initial point $(t,x(t))|_{t=0}=(0,x_{0})\in \R\times \R^{n}$. It is easy to see that $\lambda_{1}^{s}\rightarrow -\oo$ since $\lambda$ is arbitrary one in $\R$, then we have $\emptyset =\Sigma_{L}(A)=\Sigma_{NED}(A)$.

Now,  we  prove  the   theorem   for
the  nontrivial case ($\Sigma_{NED}(A)\ne \emptyset $ and $\Sigma_{NED}(A)\ne \R$).
Choosing $\gamma\in \rho_{NED}(A)$, define
\[S_{\gamma}:=\left\{(\tau,\xi)\in \R \times \R^{n}: \sup_{t\ge \tau}\{\|\Phi(t,\tau)\xi\|e^{-\gamma t}\}e^{-\varepsilon \tau}<\oo\right\},\]
and
\[U_{\gamma}:=\left\{(\tau,\xi)\in \R \times \R^{n}: \sup_{t\le \tau}\{\|\Phi(t,\tau)\xi\|e^{-\gamma t}\}e^{-\varepsilon \tau}<\oo\right\}.\]

Then for any $\gamma_{j}\in \rho_{NED}(A)$, i.e.,
\[b_{j}<\gamma_{j}<a_{j+1},\quad {\rm for}~j=1,\cdots,n-1,\]
the intersection
\[W_{j}=U_{\gamma_{j-1}}\cap S_{\gamma_{j}},\quad {\rm for}~j=1,\cdots,n-1\]
forms a  linear   integral  manifold  of \x{b1} with $\dim W_{j}\ge 1$ (see \cite{chu-15} for details). Let $ \lambda$ be an arbitrary  point  in $(b_{j}, a_{j+1})$, thus $\lambda\in \rho_{NED}(A)$, and there exist an invariant projection  $P$, constants $\alpha>0$, $M>0$, and $\varepsilon \in [0,\alpha)$ such that
 \[
\|\Phi_{\lambda}(t)P\Phi_{\lambda}^{-1}(s)\|\le M e^{-\alpha(t-s)}e^{\varepsilon s}, \quad {\rm for} \quad 0\le s \le t,
\]
or  equivalently,
\[
\|\Phi(t)P\Phi^{-1}(s)\|\le M e^{(\lambda-\alpha)(t-s)}e^{\varepsilon s}, \quad {\rm for} \quad 0\le s \le t.
\]
Set $s=0$, the inequality above implies that \[\|\Phi(t)P\Phi^{-1}(0)x_{0}\|\le M \|x_{0}\| e^{(\lambda-\alpha)t} \quad {\rm for} \quad t \ge 0
\] with any initial point $(t,x(t))|_{t=0}=(0,x_{0})\in W_{j}$. Now we can prove that $\lambda_{j}^{s}\le b_{j}$ since $\lambda$ is arbitrary one in $(b_{j}, a_{j+1})$.

Similarly, choosing $ \lambda\in (b_{j-1}, a_{j})$, thus $\lambda\in \rho_{NED}(A)$, and there exist an invariant projection  $Q=I_{n}-P$, constants $\alpha>0$, $M>0$, and $\varepsilon \in [0,\alpha)$ such that
\[\|\Phi_{\lambda}(t)Q\Phi_{\lambda}^{-1}(s)\|\le M e^{\alpha(t-s)}e^{\varepsilon s}, \quad {\rm for} \quad 0\le t \le s,\]
or  equivalently,
\[\|\Phi(t)Q\Phi^{-1}(s)\|\le M e^{(\lambda+\alpha)(t-s)}e^{\varepsilon s}, \quad {\rm for} \quad 0\le t \le s.\]
Set $t=0$, the inequality above implies that
\[\|x_{0}\|=\|\Phi(0)Q\Phi^{-1}(s)x(s)\|\le M \|x(s)\|e^{-(\lambda+\alpha-\varepsilon)s},\quad {\rm for} \quad s \ge 0\] with any initial point $(t,x(t))|_{t=0}=(0,x_{0})\in W_{j}$. Then we have
\[\|x(s)\| \ge e^{(\lambda+\alpha-\varepsilon)s}\|x_{0}\|,\quad {\rm for} \quad s \ge 0,\]
which means that $\lambda_{j}^{i}\ge \lambda+\alpha-\varepsilon\ge \lambda$ due to the fact $\varepsilon \in [0,\alpha)$.  Now we can prove that $\lambda_{j}^{i}\ge a_{j}$ since $\lambda$ is arbitrary one in $(b_{j-1}, a_{j})$.  \hspace{\stretch{1}}$\Box$
 }

 \vspace{6pt}
The next connection concentrate on the perturbation results of spectrums $\Sigma_{L}(A)$, $\Sigma_{ED}(A)$ and $\Sigma_{NED}(A)$.
It is well known that exponential dichotomy of \x{b1} remains unchanged with a  small perturbation, which is called \emph{ roughness} (see e.g., \cite[pp. 34]{cop78} for details), i.e., for a perturbed system
 \be\lb{b10}\dot{x}=(A(t)+B(t))x\ee
 with $\|B(t)\|\le \delta$ for some sufficiently small  $ \delta >0$,  the perturbed equation \x{b10} has also an exponential dichotomy. Thus, $\Sigma_{ED}(A)$ is stable under small perturbation, since the shifted system does not change the stability of exponential dichotomy. In \cite{bv-08}, Barreira and Valls show that the perturbation with the coefficient matrix being exponentially decaying, i.e.,  the linear perturbed system \x{b10}  has also a nonuniform  exponential dichotomy, while $\|B(t)\|\le \delta e ^{-\varepsilon t}$ for some sufficiently small $ \delta >0,~\varepsilon \in [0,\alpha)$.  Thus, $\Sigma_{NED}(A)$ is stable with the perturbation of the coefficient matrix being exponentially decaying.

The stability theory of Lyapunov spectrum $\Sigma_{L}(A)$ is more complicated than $\Sigma_{ED}(A)$ and $\Sigma_{NED}(A)$, we first mention that  it is not enough to ensure the stability of of Lyapunov exponents for a general system even if  for  a regular system with different Lyapunov exponents. Example from \cite[p. 171]{ad-95} shows that a two dimension system \[
\dot{x}_{1}=(1+\frac{\pi}{2}sin(\pi \sqrt{t}))x_{1},\quad \dot{x}_{2}=0
\]
has distinct Lyapunov exponents $\lambda_{1}=1$ and $\lambda_{1}=0$. However, the Lyapunov exponents of this system are not stable.

A general condition called integral separateness (see, e.g., \cite{ad-95}), which is introduced and improved by by Bylov, Vinograd, Izobov, Grobman, Million\v{s}\v{c}ikov and several others \cite{bd-07, by-54, by-65, by-66, by-69, M-69,M-69-2, vi-60}, is generally used to guarantee the stability of of Lyapunov exponents. Now we introduce the definitions of weak integral separateness, which extend the concept of integral ones.

\begin{definition} \lb{def32} {\rm(see \cite[Def. 2.2]{zc})} The continuous functions $g_{i}, i=1, \ldots, n$, are said to be {\rm weakly
integrally separated} if for $i=1, \ldots, n-1$,
there exist some costants $a, b \geq 0$ and $d\in \R$ such that
\[
\int^{t}_{s}(g_{i+1}(\tau)-g_{i}(\tau)) d \tau\geq a(t-s)-b s + d, \quad t\ge s \ge 0.
\]
\end{definition}

\begin{definition} {\rm(see \cite[Def. 2.3]{zc})} Let $\Phi(t)=(\Phi_{1}(t),\ldots,\Phi_{n}(t))$ be a fundamental matrix solution of {\rm \x{b1}}. Then, system {\rm \x{b1}} is said to be {\rm weakly integrally separated} if for $i=1, \ldots, n-1$, there exist some constants $a, b \geq 0$ and $D > 0$ such that
\be\lb{b11}
\frac{\|\Phi_{i+1}(t)\|}{\|\Phi_{i+1}(s)\|} \cdot \frac{\|\Phi_{i}(s)\|}{\|\Phi_{i}(t)\|}\geq D e^{a(t-s)-b s}, \quad t\ge s \ge 0.
\ee
\end{definition}

\x{b1} is called integrally separated if $a>0$ and $b=0$ in \x{b11} (see e.g. \cite[Definition 5.3.2]{ad-95} and \cite{dv-02}). Obviously, integral separateness implies week integral separateness due to the fact $b\ge 0$, but not vice versa. Indeed, \x{b12}  is weakly integrally separated but not integrally separated.

The following two theorems present the necessary and sufficient conditions of the stability of Lyapunov exponents, and therefore the corresponding stability of  $\Sigma_{L}(A)$.

\begin{theorem}\lb{lem1} {\rm(see \cite[Thm. 5.4.7]{ad-95} and \cite{by-69})}
Assume that the system {\rm\x{b1}} has distinct Lyapunov exponents $\lambda_{1} > \cdots > \lambda_{n}$.
Then they are stable if and only if there exists a fundamental matrix solution with  integrally separated columns.
\end{theorem}

\begin{theorem}\lb{lem2} {\rm(see \cite{zc})}
Assume that the system {\rm\x{b1}} with nonuniformly bounded growth  has distinct Lyapunov exponents $\lambda_{1} > \cdots > \lambda_{n}$.
Then they are stable with the perturbations of the coefficient matrix being exponentially decaying i.e., for a perturbed system {\rm \x{b10}}
 with $\|B(t)\|\le \delta e ^{-\varepsilon t}$ for some $ \delta >0,~\varepsilon \in [0,\alpha)$, the Lyapunov exponents of system {\rm\x{b1}} are stable
 if and only if there exists a fundamental matrix solution with weakly integrally separated columns.
\end{theorem}

 From the analysis above, we have the following perturbation results about $\Sigma_{L}(A)$, $\Sigma_{ED}(A)$ and $\Sigma_{NED}(A)$.

\begin{proposition} \lb{pro2}
   For an $n$-dimensional linear system {\rm\x{b1}}.
   \begin{enumerate}
     \item Given a sufficiently small parameter $ \delta >0$, such that $\|B(t)\|\le \delta$, then the perturbed system {\rm\x{b10}} having the following properties:
         \begin{enumerate}
           \item $\Sigma_{ED}(A)$ is stable under the perturbation $\|B(t)\|\le \delta$;
           \item {\rm\x{b1}} is integrable separated $\Leftrightarrow$ $\Sigma_{L}(A)$ is stable under the perturbation $\|B(t)\|\le \delta$;
         \end{enumerate}

     \item Given a sufficiently small parameter  $\delta >0$, and $\varepsilon \in [0,\alpha)$, such that $\|B(t)\|\le \delta e ^{-\varepsilon t}$, then the purturbed system {\rm\x{b10}} having the following properties:
         \begin{enumerate}
           \item $\Sigma_{NED}(A)$ is stable under the perturbation $\|B(t)\|\le \delta e ^{-\varepsilon t}$;
           \item {\rm\x{b1}} is weakly integrable separated $\Leftrightarrow$ $\Sigma_{L}(A)$ is stable under the perturbation $\|B(t)\|\le \delta e ^{-\varepsilon t}$;
         \end{enumerate}
   \end{enumerate}
 \end{proposition}

From the Proposition \ref{pro2}, one can find perturbation results of $\Sigma_{L}(A)$, $\Sigma_{ED}(A)$ and $\Sigma_{NED}(A)$ and some connection between these three spectrums, which can be used to assure the correctness of the numerical works in the next section.   To further explore the relationship of these three spectrums, we focus on a special case: full spectrum.
  In \cite{bs-00}, Bodine and Sacker presented that the system \x{b1} with full exponential dichotomy spectrum is integrally separated. The converse does not
hold (see \cite[pp. 193]{pa-82-2} for details). Thus the following  relationship holds:
\[\Sigma_{ED}(A) ~{\rm is ~full} \Longrightarrow ~{\rm Integral ~Separation}\]

More recently, we have proved in \cite{zc} that the system \x{b1} with full nonuniform exponential dichotomy spectrum is weakly integrally separated, and the converse can hold true if we consider some additional condition (see \cite[Theorem 1.2]{zc} for details), which can also be used to prove the converse part from integrally separated to full exponential dichotomy spectrum. This means that
\[\Sigma_{NED}(A) ~{\rm is ~full} \Longrightarrow ~{\rm Weakly~integral ~Separation}\]

Note that $\Sigma_{NED}(A)\subset \Sigma_{ED}(A)$, and integral separateness always implies week integral separateness.  Combine these  relationship with the approaches above, we will have the following
chain of implications
\[
  \begin{array}{ccc}
    \Sigma_{ED}(A) ~{\rm is ~full}  & \Rightarrow & {\rm Integral ~Separation}\\
    \Downarrow &  & \Downarrow \\
    \Sigma_{NED}(A) ~{\rm is ~full}& \Rightarrow &{\rm Weakly~integral ~Separation}
  \end{array}
\]

\begin{remark}
  If $\Sigma_{ED}(A)$ is full, the perturbation results of $\Sigma_{L}(A)$ and $\Sigma_{ED}(A)$ are the same. Similarly, if $\Sigma_{NED}(A)$ is full, the perturbation results of $\Sigma_{L}(A)$ and $\Sigma_{NED}(A)$ are the same.
\end{remark}

\section{\bf{Numerical computation of  $\Sigma_{ED}(A)$ and $\Sigma_{NED}(A)$ without bounded condition.}}
\setcounter{equation}{0} \noindent
Recall  that \x{b1} has \emph{bounded growth} (see \cite{sie-02} and \cite[pp. 8]{cop78}) if and only if there exist constants $K>0$, $\tilde{a}>0$  such that
\be\lb{d1}
\|\Phi(t)\Phi^{-1}(s)\|\le K e^{\tilde{a}|t-s|}, \quad {\rm for} \quad  t, s \ge 0.
\ee

However, the notion of bounded growth demands considerably from
the dynamics and it is of considerable interest to look for more general types of
hyperbolic behavior. We now present an en example without uniform bounded growth.

\begin{example}\lb{exm41}
  The scalar equation
\be\lb{d2}\dot{x}=t (\sin t+1) x\ee
has no  uniform bounded growth.
\end{example}
\prf{It is easy to verify that
\beaa \exp\left(\int_{s}^{t} \tau(\sin \tau +1) d \tau\right)\EQ \exp(-t\cos t  + s \cos s +\sin t -\sin s +t-s)\\
\EQ \exp(2(t-s)- t(\cos t+1)+s(\cos s+1)+(\sin t -\sin s)).\\
\LE \exp(2(t-s) +2s + 2).
\eeaa
Furthermore, if $t=2k\pi+\pi$ and $s=2k\pi$ with $k\in \N$, then
\be\lb{d3}\exp\left(\int_{s}^{t} \tau(\sin \tau+1) d \tau\right) =\exp(2(t-s) +2s).\ee
Thus \x{d2}  has no uniform bounded growth due to the fact that the perturbation $2s$ in \x{d3} can not be eliminated. This means that bounded growth \x{d1} is not satisfied. \hspace{\stretch{1}}$\Box$
}

\begin{remark}\lb{remfj1} The numerical method proposed in \cite{dv-02,dv2-02} is not quite right without bounded condition.
In fact, the computational procedure to approximate dichotomy spectrum is based on the separateness of Steklov function, which is equivalent to integral separateness under the condition of bounded (See Lemma 5.4.1 in \cite{ad-95}).

For example, consider
the following scalar equation
\be \lb{a1} \dot{x}=a(t)x
\ee
with $a(t)$ is continuous. The computational procedure to approximate dichotomy spectrum of \x{a1} is as follows. Given $H>0$, and $T>t_{0}>0$.  Let $a_{H}=\frac{1}{H}\int_{t}^{t+H}a(\tau) d \tau$ for $t\in [t_{0},T]$. One can compute
\[ \overline{a}=\sup_{t_{0}\le t \le T-H} a_{H} \quad {\rm and} \quad
\underline{a}=\inf_{t_{0}\le t \le T-H} a_{H}
\]
and use $[\underline{a}, \overline{a}]$ as an approximation to dichotomy spectrum of \x{a1}.

Let $a(t)=t (\sin t+1)$ as in \x{d2}. Let $t_{1}=2k_{1}\pi$, and $H_{1}=2k_{2}\pi+\pi$  with $k_{1}, k_{2}\in \N$,  it follows easily from \x{d3} that
\[\frac{1}{H_{1}}\int_{t_{1}}^{t_{1}+H_{1}} \tau (\sin \tau+1) d \tau= 2+ \frac{t_{1}}{H_{1}},\]
Similarly, let $t_{1}=2k_{1}\pi+\pi$, and $H_{1}=2k_{2}\pi$  with $k_{1}, k_{2}\in \N$, we have
\[\frac{1}{H_{2}}\int_{t_{2}}^{t_{2}+H_{2}} \tau (\sin \tau+1) d \tau=- \frac{t_{2}}{H_{2}},\]
Hence, the dichotomy spectrum of \x{d2} is $\R$, due to the fact that $a(t)=t (\sin t+1)$  is not bounded.
\end{remark}

In order to present the numerical computation of spectral intervals $\Sigma_{ED}(A)$ and $\Sigma_{NED}(A)$, we need to introduce the following definition to extend the known results of bounded growth to nonuniformly bounded growth.

\begin{definition} \lb{def41}
{\rm(see~\cite[Def. 2.9]{chu-15})}
  We say that {\rm\x{b1}}  has a {\rm nonuniformly bounded growth} if there exist constants $K>0$, $\tilde{a}>0$ and $\tilde{b}\ge 0$ such that
\begin{equation*}
\|\Phi(t)\Phi^{-1}(s)\|\le K e^{\tilde{a}|t-s|}e^{\tilde{b} s}, \quad {\rm for} \quad t,s \ge 0,
\end{equation*}
where $\Phi(t)$ is a fundamental matrix of {\rm\x{b1}} .
\end{definition}

Recall that the function
\be\lb{d5} f^{H}(t)=\frac{1}{H}\int_{t}^{t+H}f(\tau) d \tau
\ee
is defined as \emph{ Steklov function} or \emph{Steklov average} (see \cite[Def. 5.4.1]{ad-95}) with step $H>0$.   Inspired by the result
\[\int_{s}^{t} \tau(\sin \tau+1) d \tau \le 2(t-s) +2s + 2\]
 in
Example \ref{exm41}, the following lemma, unlike the work in \cite[Lemma 5.4.1]{ad-95}, is to  investigate the necessary and sufficient condition of weak integral separateness under the condition of nonuniform bounded growth, i.e., \be\lb{d8}\int_{s}^{t}|f(\tau)| d \tau \le \tilde{a}|t-s|+ \tilde{b} s +\tilde{d} \quad t, s\ge 0,\ee with $\tilde{a}, \tilde{b}>0$ and $\tilde{d}\in \R$, since the fundamental matrix solution of $\dot{x}=f(t)x$ satisfies
\[ |\Phi(t)\Phi^{-1}(s)|\le e^{\int_{s}^{t}|f(\tau)| d \tau}\le e^{\tilde{d}} e^{\tilde{a}|t-s|+ \tilde{b} s }.\]

\begin{lemma} \lb{lem41}
  Assume that $f_{1}(t)$ and $f_{2}(t)$ are nonuniformly bounded growth functions, i.e., {\rm \x{d8}} holds true. Then  $f_{1}(t)$ and $f_{2}(t)$ are weakly integrally separated with $a, b>0$ if and only if for sufficiently large $H >> t$, the Steklov functions are separated in the standard sence:
        \be \lb{d6} f_{2}^{H}(t)-f_{1}^{H}(t)\ge M\ee
        for some constant $M>2\tilde{b}>0$;
  \end{lemma}
\prf{From the proof of Lemma 5.4.1 in \cite{ad-95}, we know that the equality
\be \lb{d9}\int_{s}^{t}f^{H}(\tau) d \tau =\int_{s}^{t}f(\tau) d \tau +I(t)-I(s)\ee
holds with
\[I(t)=\frac{1}{H}\int_{t}^{t+H}f(y)d y \int_{y-H}^{t}d x.\]

 Using \x{d8}, we have
\[|I(t)|\le \frac{1}{H}\int_{t}^{t+H}(t-y+H)|f(y)|d y \le \tilde{a}H+ \tilde{b} t +\tilde{d}.\]
Thus it follows from \x{d6} and \x{d9} that
\beaa\int_{s}^{t}(f_{2}(\tau)-f_{1}(\tau)) d \tau \EQ \int_{s}^{t}(f^{H}_{2}(\tau)-f^{H}_{1}(\tau)) d \tau  -I_{2}(t)+I_{2}(s) +I_{1}(t)-I_{1}(s)\\
\GE (M-2\tilde{b})(t-s)-4\tilde{b}s-4(\tilde{a}H+\tilde{d}),
\eeaa
this implies that the functions $f_{1}(t)$ and $f_{2}(t)$ are weakly integrally separated.

Conversely, assume that  the functions $f_{1}(t)$ and $f_{2}(t)$ are weakly integrally separated, then we have
\[
\int^{t}_{s}(f_{2}(\tau)-f_{1}(\tau)) d \tau\geq a(t-s)-b s + D, \quad t\ge s \ge 0
\]
for $a, b>0$ and $D\in \R$. Thus the difference of Steklov functions is
\be\lb{zj1}
f^{H}_{2}(t)-f^{H}_{1}(t)=\frac{1}{H}\int^{t+H}_{t}(f_{2}(\tau)-f_{1}(\tau)) d \tau\geq a-b \frac{t}{H} + \frac{D}{H}.
\ee
Hence \x{d6} holds with $H >> t$. \hspace{\stretch{1}}$\Box$
}

\begin{remark}\lb{rem41}
 From the proof of Lemma \ref{lem41},  $H$ in Steklov function \x{d5} must be chosen such that  $H >> t$, or else, $\frac{t}{H}$ in \x{zj1} can not be  ignored, which is completely different from those in \cite[Lemma 5.4.1]{ad-95}.

Moreover, in the actual calculation process, we need the condition $t >> H$ to find the the effect of the nonuniform item. In fact, it follows from Definition \ref{def32} that
\[\frac{H}{t}\left|f_{2}^{H}(t)-f_{1}^{H}(t)\right| = \frac{1}{t}\left|\int^{t+H}_{t}(f_{2}(\tau)-f_{1}(\tau)) d \tau\right|
\ge b-a\frac{H}{t}-\frac{|D|}{t} \ge 0
\]  with $t >> H$. This means that for sufficiently large $t >> H$, the Steklov functions satisfy the inequality
    \be \lb{d7}\frac{H}{t}\left|f_{2}^{H}(t)-f_{1}^{H}(t)\right|\ge N\ee
        for some constant $N>0$ if the nonuniform term does exist, and this effect does not appear in \cite[Lemma 5.4.1]{ad-95} with $f_{1}(t)$ and $f_{2}(t)$ are  bounded, or even if $f_{1}(t)$ and $f_{2}(t)$ are uniformly bounded growth functions.
\end{remark}

 In this paper, we always assume that \x{b1} is weakly integrally separated.  Note that a weakly integrally separated system is invariant under Lyapunov transformation, and a weakly integrally separated system is kinematically similar to a diagonal one by using the Lyapunov transformation (see \cite{zc}).  So for a diagonal system, or for any system which can be reduced to a diagonal system through a Lyapunov transformation, our approach for approximating $\Sigma_{L}(A)$ of \x{b1} under the condition of nonuniform bounded growth is the same as $A$ is bounded in \cite{dv-02}.  Hence, on a finite time interval, our computational procedure for $\Sigma_{L}(A)$ is as follows. Consider a diagonal system
  \be \lb {d10}\dot{y}=diag [a_{1}(t),\ldots, a_{n}(t)]y.\ee
 Given  constants $T_{1}$, $T_{2}>0$, such that $t\in [T_{1},T_{2}]$ with $T_{2}>>T_{1}>0$. Let $\lambda_{j}(t)=\frac{1}{t}\int_{0}^{t}a_{j}(\tau)d\tau$, and compute
\[ \overline{\lambda}_{j}=\sup_{T_{1}\le t \le T_{2}} \lambda^{t}_{j} \quad {\rm and} \quad
\underline{\lambda}_{j}=\inf_{T_{1}\le t \le T_{2}} \lambda^{t}_{j}
\]
and use $[\underline{\lambda}_{j}, \overline{\lambda}_{j}]$ as an approximation to $[\lambda_{j}^{i},\lambda_{j}^{s}]$.

However, under the condition of nonuniform bounded growth, the procedure for approximation of $\Sigma_{NED}(A)$  is essentially different from the approximation of $\Sigma_{ED}(A)$ in \cite{dv-02}. Indeed, the nonuniform item can cause catastrophic failure in the computation when the approximation scheme in \cite{dv-02} is applied here since the nonuniform item can not be eliminated  (see Remark \ref{remfj1} for details).

\begin{definition}
   the weak integral separation spectrum is
   \[\Sigma_{WIS}=\bigcup_{j=1}^{n}\Lambda_{j},\] where $\Lambda_{j}=\Lambda_{j}^{+}\bigcap\Lambda_{j}^{-}$ is a closed interval, with $\Lambda_{j}^{+}$
   corresponding to the $j$th diagonal planar systems
  \be \lb{d11} \dot{y}_{j}=\left(
                    \begin{array}{cc}
                      \lambda & 0 \\
                      0 & a_{j}(t) \\
                    \end{array}
                  \right)y_{j},
  \ee
 and $\Lambda_{j}^{-}$
   corresponding to the $j$th diagonal planar systems
  \be \lb{d12} \dot{y}_{j}=\left(
                    \begin{array}{cc}
                     a_{j}(t) & 0 \\
                      0 & \lambda  \\
                    \end{array}
                  \right)y_{j}
  \ee
for each $j=1,\cdots, n$, which
   are given by
 \[\Lambda_{j}^{+}=\{\lambda \in \R: {\rm~\x{d11} ~is ~not~ weakly ~integrally ~separated}\}.
   \]
   and
    \[\Lambda_{j}^{-}=\{\lambda \in \R: {\rm~\x{d12} ~is ~not~ weakly ~integrally ~separated}\}.
   \]
\end{definition}

The following theorem mimics the classical one about integral separation spectrum \cite[Theorem 2.29]{dv2-02}, but for the variation of spectrum has been extended from uniform item to the nonuniform ones.

\begin{theorem}\lb{thm41}
  For {\rm \x{d10}}, $\Sigma_{WIS}=\Sigma_{NED}$.
\end{theorem}

\prf { Given $\lambda \in \R$, if $\lambda \notin \Sigma_{NED}$. It follows from \x{b8}-\x{b9} that there  exist constants $\alpha>0$, $M>0$, and $\varepsilon \in [0,\alpha)$ such that either
\be \lb{d13} e^{\int_{s}^{t}a_{j}(\tau)d\tau}e^{-\lambda(t-s)} \le M e^{-\alpha(t-s)}e^{\varepsilon s}, \quad {\rm for} \quad 0\le s \le t,
\ee
or
\be \lb{d14} e^{\int_{s}^{t}a_{j}(\tau)d\tau}e^{-\lambda(t-s)} \le M e^{\alpha(t-s)}e^{\varepsilon s}, \quad {\rm for} \quad 0\le t \le s.
\ee
If \x{d13} holds then \x{d12} is weakly integrally separated. If \x{d14} holds then
\x{d11}  is weakly integrally separated. This means that $\lambda \notin \Sigma_{NED}\Rightarrow \lambda \notin \Sigma_{WIS}$. Conversely, if $\lambda \notin \Sigma_{WIS}$, then for all $j=1,\cdots, n$, either \x{d11} or \x{d12} is weakly integrally separated and hence either  \x{d13} or \x{d14} hold. \hspace{\stretch{1}}$\Box$
}

\begin{theorem}\lb{thm42}
  Consider the diagonal system {\rm\x{d10}} with nonuniformly bounded growth, i.e., for $j=1,\ldots,n$,
  \be\lb{d15}\int_{s}^{t}|a_{j}(\tau)|d \tau \le \tilde{a}(t-s)+ \tilde{b} s +\tilde{d} \quad t\ge s\ge 0,\ee with $\tilde{a}, \tilde{b}>0$ and $\tilde{d}\in \R$.
  Let
  \[ \alpha_{j}^{H} =\inf_{t}\frac{1}{H}\int_{t}^{t+H}a_{j}(\tau) d \tau
  \quad {\rm and}\quad \beta_{j}^{H} =\sup_{t}\frac{1}{H}\int_{t}^{t+H}a_{j}(\tau) d \tau
  \]
 with any given $H >0$.
  Then, for each $j=1,\ldots,n$, $\Lambda_{j}\subseteq [\alpha_{j}^{H}, \beta_{j}^{H}]$.  Moreover, assume that {\rm \x{d11}} and {\rm \x{d12}} are weakly integrally separated respectively with $a,~b>0$. Then for $H>>t$ sufficiently large,  $[\alpha_{j}^{H},\beta_{j}^{H}]\subseteq \Lambda_{j}$ for $j=1,\ldots,n$.
\end{theorem}

\prf{ First, let $\lambda > \beta_{j}^{H}$. Thus, there exists $M_{j}>0$ such that
\be\lb{d16}
\int_{t}^{t+H}(\lambda-a_{j}(\tau))d \tau \ge M_{j} H\quad \forall ~t
\ee
In order to prove that $\lambda$ and $a_{j}$ are weakly integrally separated, it suffices to present that there exist some costants $a, b > 0$ and $D\in \R$ such that
\be \lb{d17}
\int^{t}_{s}(\lambda-a_{j}(\tau)) d \tau\geq a(t-s)-b s - D, \quad t\ge s \ge 0.
\ee
We will verify \x{d17} with $a=M_{j}$, $b=\tilde{b}$, $D=H(\lambda +\tilde{a})+\tilde{b} H+ \tilde{d}$. In fact, it is easy to see that \x{d17} holds with $a=M_{j}$ for all $t$ and $s$ with $t=s+H$ because of the inequality \x{d16}. Now consider the case with $t<s+H$, we can rewrite the left hand side of \x{d17} as
\[\int^{t}_{s}(\lambda-a_{j}(\tau)) d \tau= \int^{s+H}_{s}(\lambda-a_{j}(\tau)) d \tau
-\int^{s+H}_{t}(\lambda-a_{j}(\tau)) d \tau,\]
and thus it follows from $s \le t<s+H$ and  \x{d15} that
\beaa\int^{s+H}_{t}(\lambda-a_{j}(\tau)) d \tau \LE (\lambda +\tilde{a})(s+H-t)+ \tilde{b} (s+H) +\tilde{d}\\
\LE bs+D,
\eeaa
and hence,
\[\int^{t}_{s}(\lambda-a_{j}(\tau)) d \tau \ge aH-bs-D \ge a(t-s)-bs-D.\]

Next, let $t>s+H$. Then, for some integer $k>1$, $t=s+kH+\rho$ with $\rho\in [0,H)$. Then we have
\beaa\int^{t}_{s}(\lambda-a_{j}(\tau))d \tau\EQ\sum_{j=0}^{k}\int_{t-(j+1)H}^{t-jH}(\lambda-a_{j}(\tau))d \tau
-\int_{t-(k+1)H}^{s}(\lambda-a_{j}(\tau))d \tau\\
\GE M_{j}(k+1)H -((\lambda +\tilde{a})(H-\rho)+\tilde{d}) - \tilde{b} s \\
\GE  a(t-s)-bs-D.
\eeaa
Therefore, $\lambda$ and $a_{j}$ are weakly integrally separated, and this means that
$\lambda \notin \Lambda_{j}$. Similarly, we can prove that $\lambda \notin \Lambda_{j}$ for any $\lambda < \alpha_{j}^{H}$, and so $\Lambda_{j}\subseteq [\alpha_{j}^{H},\beta_{j}^{H}]$ for any given $H>0$.

 Conversely, assume that $\lambda \notin \Lambda_{j}$, then  $\lambda$ and $a_{j}$ or $a_{j}$ and $\lambda$ are weakly integrally separated. Suppose that $\lambda$ and $a_{j}$ are weakly integrally separated with constants $a, b>0$ and $D\in \R$. Thus for any given $t\in \R$, choosing $H>>t$ such that
 \be \lb{d19}
 \frac{1}{H}\int_{t}^{t+H}(\lambda-a_{j}(\tau))d \tau  \ge a- b\frac{t}{H}-\frac{D}{H}> \frac{a}{2}>0,
 \ee
 and so $\lambda> \beta_{j}^{H}$. Similarly, we can prove that $\lambda<\alpha_{j}^{H}$.  Therefore,  for $H>>t$ sufficiently large,  $[\alpha_{j}^{t},\beta_{j}^{t}]\subseteq \Lambda_{j}$ for $j=1,\ldots,n$.
\hspace{\stretch{1}}$\Box$
}

  Unlike the result of  \cite{dv-02}, $[\alpha_{j}^{H},\beta_{j}^{H}]\subseteq \Lambda_{j}$ does not hold in general for the nonuniform bounded case. In fact, it follows from \x{d19} that $\alpha_{j}^{H}$ and/or $\beta_{j}^{H}$ can be unbounded if $t>> H$ (see Remark \ref{remfj1} for details).

Now, to obtain a computational procedure on a finite time interval for $\Sigma_{NED}$ out of Theorem \ref{thm42}, we need to verify whether the condition \x{d7} holds or not. If \x{d7} holds with $N>0$, then the nonuniform term in  \x{d10}  does exist, and this means that $\Sigma_{ED}=\R$ since the nonuniform term can not be  eliminated (see Example 2.1 in \cite{chu-15} for details). Otherwise, $N=0$ in \x{d7} means that there is no nonuniform term in  \x{d10}, that is, $\Sigma_{NED}=\Sigma_{ED}$.

Hence, on a finite time interval, our computational procedure is as follows. First, following the ideas in Remark \ref{rem41}, we compute the size of bias of nonuniform item $b$ in \x{d17}. Given any $H>0$, there exist constants $T_{1}$, $T_{2}>0$, such that $t\in [T_{1},T_{2}]$ and $T_{1}>>H$. Let $b^{t}_{j}=\frac{1}{t}\int_{t}^{t+H}a_{j}(\tau) d \tau$ for $t\in [T_{1},T_{2}]$. Then we can compute
\be\lb{b21} \overline{b}_{j}=\sup_{T_{1}\le t \le T_{2}} |b^{t}_{j}|,
\ee
and use $\overline{b}_{j}$  to represent the bias of the nonuniform item $b$ in \x{d17}.
if $0<\overline{b}_{j}<\epsilon<<1$ for some $\epsilon>0$,  there is no nonuniform term in   \x{d10}, or else, the nonuniform term in  \x{d10}  does exist.

Now we compute the spectrum $\Sigma_{NED}$ if the nonuniform item $b$ in \x{d17} far away from zero, otherwise, we can follow the idea in \cite{dv-02} to compute the spectrum $\Sigma_{ED}$ since there is no nonuniform item $b$ in \x{d17}. Thus, it follows from Theorem \ref{thm42} that a computational procedure for $\Sigma_{NED}$ on a finite time interval is as follows. Given  constants $T_{1}$, $T_{2}>0$, and $H>0$, such that $t\in [T_{1},T_{2}]$ and $H>>T_{2}$. Let $a^{t}_{j}=\frac{1}{H}\int_{t}^{t+H}a_{j}(\tau) d \tau$ for $t\in [T_{1},T_{2}]$. Then we compute
\be\lb{b20} \overline{a}_{j}=\sup_{T_{1}\le t \le T_{2}} a^{t}_{j} \quad {\rm and} \quad
\underline{a}_{j}=\inf_{T_{1}\le t \le T_{2}} a^{t}_{j}
\ee
and use $[\underline{a}_{j}, \overline{a}_{j}]$ as an approximation to $[\alpha_{j}^{H}, \beta_{j}^{H}]$.

%
%

\section{\bf{Example and Numerical Simulation}}
\setcounter{equation}{0} \noindent
In this Section, we  consider a planar problem, which satisfies the condition of nonuniform  bounded growth. In this case, we approximate the spectral intervals of $\Sigma_{L}(A)$, $\Sigma_{ED}(A)$ and $\Sigma_{NED}(A)$ and compute the bias of nonuniform item.

\begin{example}
  Consider a planar system
  \be\lb{e1} \left(\begin{array}{c}
         \dot{x}_{1} \\
         \dot{x}_{2}
       \end{array}\right)=\left(
  \begin{array}{cc}
    sin (ln(t)) + cos (ln(t)) & 0\\
    0 & \omega_{1}-\omega_{2}t\sin t \\
  \end{array}
\right)\left(\begin{array}{c}
         x_{1} \\
         x_{2}
       \end{array}\right)
\ee
with $\omega_{1}>\omega_{2}>0$.
The problem is designed so that the coefficient matrix of \x{e1} is nonuniform bounded growth. Note that the solution of \x{e1} is
\[ \left(\begin{array}{c}
         x_{1} \\
         x_{2}
       \end{array}\right)=\left(\begin{array}{c}
         \exp(t sin (ln(t))-t_{0} sin (ln(t_{0}))) \cdot x_{1}(t_{0}) \\
         \exp(\omega_{1}(t-t_{0})+\omega_{2}t\cos t  -\omega_{2} t_{0} \cos t_{0}- \omega_{2}\sin t+ \omega_{2}\sin t_{0}) \cdot  x_{2}(t_{0})
       \end{array}\right).
\]
Hence, $\Sigma_{L}=[-1,1]\cup [\omega_{1}-\omega_{2},\omega_{1}+\omega_{2}]$.
Moreover, it follows from Example 6.2 of \cite{dv-02} and Example 2.1 of \cite{chu-15}
that  $\Sigma_{ED}=[-\sqrt{2},\sqrt{2}] \cup \R$ and $\Sigma_{NED}=[-\sqrt{2},\sqrt{2}] \cup [\omega_{1}-\omega_{2},\omega_{1}+\omega_{2}]$.
 In the actual computation, here we choose $\omega_{1}=4$, and $\omega_{2}=2$ for \x{e1}.

\end{example}

Our numerical results of of $\Sigma_{L}$ are listed in Table 1. In this table we specify the values of $T_{1}$ and $T_{2}$, and calculate the approximations of the two  spectral intervals of  $\Sigma_{L}$. The results in Table 1 show that the approximations are quite accurate for these three time intervals $[T_{1}, T_{2}]$.

\begin{table}[htbp]
\tabcolsep 0pt \caption{Approximation of $\Sigma_{L}$} \vspace*{-12pt}
\begin{center}
\def\temptablewidth{0.8\textwidth}
{\rule{\temptablewidth}{1pt}}
\begin{tabular*}{\temptablewidth}{@{\extracolsep{\fill}}cccc}
$T_{1}$ & $T_{2}$ & $[\underline{\lambda}_{1}, \overline{\lambda}_{1}]$ &$[\underline{\lambda}_{2}, \overline{\lambda}_{2}]$  \\   \hline
              $1.E2$    & $1.E4$        &     $[-1.0098, 1.0004]$~~   & $[2.0019,6.0000]$
\\   \hline
              $1.E2$      & $1.E6$        &     $[-1.0060, 1.0004]$   & $[2.0000,6.0000]$
\\   \hline
              $1.E4$     & $1.E6$         &     $[-1.0000, 0.9487]$   & $[2.0000,6.0000]$
       \end{tabular*}
       {\rule{\temptablewidth}{1pt}}
       \end{center}
       \end{table}

\begin{table}[htbp]
\tabcolsep 0pt \caption{Bias of Nonuniform item} \vspace*{-12pt}
\begin{center}
\def\temptablewidth{0.8\textwidth}
{\rule{\temptablewidth}{1pt}}
\begin{tabular*}{\temptablewidth}{@{\extracolsep{\fill}}ccccc}
$H$ & $T_{1}$ & $T_{2}$ & $\overline{b}_{1}$ & $\overline{b}_{2}$  \\   \hline
     $1.E2$    &     $1.E4$    & $1.E5$        &     $0.0013$ & $1.0949$\\   \hline
     $1.E3$ &         $1.E6$      & $1.E7$           & $1.2284\times 10^{-4}$ & $1.8760$
\\   \hline
      $1.E4$  &      $1.E6$     & $1.E7$    & $1.2882\times 10^{-5}$       & $ 1.0500$
       \end{tabular*}
       {\rule{\temptablewidth}{1pt}}
       \end{center}
       \end{table}

\begin{table}[htbp]
\tabcolsep 0pt \caption{Approximation of $\Sigma_{ED}$} \vspace*{-12pt}
\begin{center}
\def\temptablewidth{0.8\textwidth}
{\rule{\temptablewidth}{1pt}}
\begin{tabular*}{\temptablewidth}{@{\extracolsep{\fill}}ccccc}
$H$ & $T_{1}$ & $T_{2}$ & $[\underline{a}_{1}, \overline{a}_{1}]$ &$[\underline{a}_{2}, \overline{a}_{2}]$  \\   \hline
     $1.E3$    &     $1.E6$    & $1.E8$        &     $[-1.4142, 1.2645]$~~   & $[-1.8707\times 10^{5}, 1.8707\times 10^{5}]$  \\   \hline
     $1.E5$ &         $1.E6$      & $1.E8$        &     $[-1.4142, 1.2323]$   & $[-3.9964\times 10^{3}, 4.0044\times 10^{3}]$
\\   \hline
      $1.E4$  &      $1.E5$     & $1.E8$         &     $[-1.4142, 1.4142]$   & $[-3.9510\times 10^{4}, 3.9511\times 10^{4}]$
       \end{tabular*}
       {\rule{\temptablewidth}{1pt}}
       \end{center}
       \end{table}

       \begin{table}[htbp]
\tabcolsep 0pt \caption{Approximation of $\Sigma_{NED}$} \vspace*{-12pt}
\begin{center}
\def\temptablewidth{0.8\textwidth}
{\rule{\temptablewidth}{1pt}}
\begin{tabular*}{\temptablewidth}{@{\extracolsep{\fill}}cccc}
$H$ & $T_{1}$ & $T_{2}$  &$[\underline{a}_{2}, \overline{a}_{2}]$  \\   \hline
     $1.E4$    &     $1.E2$    & $1.E3$          & $[1.6649, 6.3694]$  \\   \hline
     $1.E6$ &         $1.E2$      & $1.E3$          & $[1.9999,5.9985]$
\\   \hline
      $1.E8$  &      $1.E3$     & $1.E4$          & $[1.9998,5.9986]$
       \end{tabular*}
       {\rule{\temptablewidth}{1pt}}
       \end{center}
       \end{table}

In Table 2 we use the computational procedure outlined in Section 4, and report on numerical results which calculate the bias of nonuniform item based on \x{b21}. It can be seen that the nonuniform items are sufficiently small for the first equation of \x{e1}  and far away from zero in the second equation of \x{e1}, which means that the first equation of \x{e1} admits an uniform exponential dichotomy, while the second equation of \x{e1} admits a nonuniform ones. Then, by calculating the the approximations of the two  spectral intervals of  $\Sigma_{ED}$ in Table 3, we can find that the second spectral interval of  $\Sigma_{ED}$  is large enough while  time intervals $[T_{1}, T_{2}]$ tends to infinity, which agrees with the theoretical result. At last, we just present the  second spectral interval of  $\Sigma_{NED}$ based upon \x{b20}  in Table 4, since the first one does not have nonuniform one. The results in this table shows that spectral interval $\Sigma_{NED}$ can be approximated accurately by letting $H>0$ large enough such that $H>>T_{2}$.

\end{document}